%% file: vertex_set_capitanio_diatta07.tex
\documentclass[11pt]{article}

\usepackage[dvips]{graphics}
\usepackage[dvips]{graphicx}
\usepackage{amsfonts}
\usepackage{amssymb}
\usepackage{amsmath}
\usepackage{amstext}
\usepackage{amsbsy}
\usepackage{amsopn}
\usepackage{amsthm}
\usepackage{amscd}
\usepackage{amsxtra}
\usepackage{upref}

\newcommand{\pp}{{\bf p}}
\newcommand{\C}{\mathbb{C}}
\newcommand{\N}{\mathbb{N}}

\newcommand{\R}{\mathbb{R}}
\newcommand{\Z}{\mathbb{Z}}

\mathchardef\varepsilon="010F
\mathchardef\epsilon="0122
\mathchardef\vartheta="0112
\mathchardef\theta="0123
\mathchardef\varrho="011A
\mathchardef\rho="0125
\mathchardef\varphi="011E
\mathchardef\phi="0127
\renewcommand
\emptyset
\varnothing

\newtheorem*{theorem*}{\bf Theorem}
\newtheorem{theorem}{\bf Theorem}
\newtheorem{lemma}{\bf Lemma}
\newtheorem*{lemma*}{\bf Lemma}
\newtheorem{proposition}{\bf Proposition}
\newtheorem*{proposition*}{\bf Proposition}

\newtheorem*{corollary*}{\bf Corollary}
\theoremstyle{definition}
\newtheorem*{definition*}{\bf Definition}
\newtheorem*{definitions*}{\bf Definitions}
\newtheorem{definition}{\bf Definition}

\newtheorem{example}{\bf Example}
\newtheorem*{example*}{\bf Example}
\theoremstyle{remark}
\newtheorem*{remark*}{\bf Remark}
\newtheorem{remark}{\bf Remark}
\newtheorem*{acknowledgements}{\bf Acknowledgements}

\author{Gianmarco Capitanio and Andr\'e Diatta}

\title{Perestroikas of vertex sets at umbilic points}

\begin{document}

\maketitle
\begin{abstract}
Mark all vertices on a curve evolving under a family of curves
obtained by intersecting a smooth surface $M$ with the 1-parameter
family of planes parallel to the tangent plane of $M$ at a point
\pp. Those vertices trace out a set, called the vertex set through
\pp. We take \pp\ to be a generic umbilic point on $M$ and describe
the perestroikas of the vertex set under generic $n$-parameter small
deformations of the surface.

Beyond the Mathematical interest on vertices of families of curves,
this work was primarily motivated by the medial representation of shapes
in Computer Vision and Image Analysis, where the behaviour of vertices
plays a crucial role in the qualitative changes of the skeleton or
Blum medial axis of curves.
\end{abstract}

%*****************************************
\section{Introduction}\label{chap:intro}
%*****************************************

Curves on surfaces play an important role in applications such as
Computer Vision, Shape Analysis, etc.
In their former work, P.J. Giblin and second
author have proposed to represent image information as a collection of
medial representations\footnote{The symmetry set of a curve
  (resp. surface) $S$ is the closure of the loci of centers of all
  circles (spheres) which are tangent to $S$ at more than one place.
  The medial axis is obtained when we only consider the circles (resp. sphere) whose
  radii equal the distance from their centers to the curve (surface).   In
  medial representation, one studies the properties and structure of
  medial axes and seeks to get information on curves (surfaces, shapes) from them.}
for the level sets of intensity (isophote curves).
They have investigated the geometry of a class of parameter families
of curves arising as a generalisation of isophote curves on surfaces.
These curves are obtained as sections of a surface by
a continuous family of planes parallel to and near the tangent
plane of the surface at a point \pp.
Such families contain singular members corresponding to the tangent
planes themselves.
Hence standard results from Singularity Theory, as in \cite{bg}, do
not apply to them.
In \cite{dg}, near any of the (elliptic, umbilic, hyperbolic,
parabolic, cusp of Gauss) points \pp\ of a generic smooth surface in
3-space, they carry out an extensive study of the behaviour of
vertices and inflexions for curves evolving near the singular member
(the curve through \pp), as well as the limits of curvatures at
vertices as the curves collapse to the singular one.
They also classify all possible arrangements
of the branches of the {\it vertex set} (set of all local patterns of
vertices) and the {\it inflexion set}.

This has been motivated, on the one hand,
by the fact that, for a curve, centers of circles of curvature at
vertices are endpoints of the so-called {\it symmetry set} of the
curve and inflexions correspond to where
a local branch of the {\it symmetry set} recedes to infinity.
Hence, from the way vertices and inflexions behave,
one can deduce a great deal of information about the local number
of branches of the symmetry set and their qualitative changes, as the
isophotes evolve.
The qualitative topological changes (or transitions) on symmetry sets
carry the information about the so-called medial axis (or skeleton)
used in medial representations of shapes.

On the other hand, the study of vertices of curves has raised up a
great interest in particular in Geometry and Singularity Theory, in
line with several problems such as the $4$-vertex Theorem, the local
geometry of surfaces and Geometry of Caustics.
These classical subjects have received a new impulsion due to the
development of Symplectic and Contact Geometry, especially in the
works of V. Arnold on Lagrangian and Legendrian Collapse and
Legendrian Sturm Theory (see \cite{arnoldlc}, \cite{arnoldams} and
\cite{arnoldfap}; see also \cite{kaz}), showing the relation between
vertices of plane curves and Lagrangian and Legendrian singularities
and Sturm-Hurwitz Theorem (see \cite{arnoldams}, \cite{ot}).
In this context, Uribe-Vargas has devoted several of his work to the
same subject.
Namely, proving a conjecture of V. Arnold, he showed that the surface
of changing four vertices for six in general $2$-parameter families of
level curves  near an umbilic point is a hypocycloidal cup
(see Problem 1993-3 of \cite{arnoldp} and \cite{uribe}).
Uribe-Vargas also proved that the bifurcation diagram has at most one
modulus.

In this paper, we describe the perestroikas of {\it the vertex sets}
at generic umbilic points of surfaces undergoing an evolution in an
$n$-parameter family of surfaces. The corresponding discriminants of
the vertex set are also studied. One of the main consequences of our
results is that, in some sense (see Theorem \ref{thm:2}),
generically the study of the discriminants of small
$n-$deformations, with $n\ge 2$, simplifies to that of
$2-$deformations. The case of generic $1$-parameter families of
surfaces is also considered in details, we draw the singular surface
of the deforming vertex set in $(2+1)$-space where the parameter is
an additionnal variable.

 As defined above, the vertex set of a surface at a point
$\pp$ is the locus of the curvature extrema on sections of the
surface by planes parallel to the tangent one at $\pp$, and these
extrema are euclidean invariants. In the conclusion (Section
\ref{conclusion}) we discuss very briefly the projective
differential geometry version of this study, as well as its possible
extension to higher degeneration.

\begin{acknowledgements}
The authors would like to thank Prof. P.J. Giblin and Dr. R.
Uribe-Vargas for very helpful discussions and comments. They are
gratefully indebted to Prof V.I. Arnold  for raising very
interesting remarks, comments and very important open questions (see
Conclusion in Section \ref{conclusion}).
\newline\noindent
The first author was supported by INdAM. A part of this work was
done while the second author was supported by the IST Programme of
the European Union (IST-2001-35443).
\end{acknowledgements}

%*****************************************
\section{Presentation of main results}\label{chap:results}
%*****************************************

Consider a smooth (embedded) curve in the Euclidian plane $\R^2$.
The {\it osculating circle} of the curve at a point \pp\
is the circle tangent to the curve at \pp\ with order at least $2$,
that is, the circle passing through $3$ infinitesimally  close points
of the curve.
The curvature of the curve at this point is the inverse of the radius
of the osculating circle.
A point of the curve is a {\em vertex} if its osculating
circle has a tangency of order at least $3$.
A vertex of a curve is a critical point of its curvature.
A vertex is {\em non-degenerate} if the corresponding osculating
circle has tangency of order $3$ with the curve.
The vertex is {\em $n$-degenerate} ($n\in\N\cup\{\infty\}$) if the
curvature of the curve has there an $A_{n+1}$ critical point.

Let us consider a smooth surface $M$ in the Euclidean $3$-space
$\R^3=\{x,y,z\}$ and a point \pp\ in $M$ considered as the origin of
$\R^3$.
Without loss of generality, we assume that the tangent plane
of the surface at \pp\ is $z=0$.
Then the surface is, at least locally, the graph $z=f_0(x,y)$ of a
smooth function $f_0:\R^2\longrightarrow \R$.

In this setting, the intersection of $M$ with planes parallel to the
tangent plane at \pp, are the level curves $f_0=k$,  where $k=0$
corresponds to the tangent plane itself. The patterns of vertices of
$f=k$, when $k$ varies, trace out a set called the {\it vertex set}
through \pp (see \cite{dg}).

We summarise some of the results in \cite{dg} as follows.

\begin{theorem}\label{thm:0}\cite{dg}
The germ of the vertex set at a point \pp\ of a generic smooth surface
in $3$-space has:
\newline\noindent
$\bullet$ four smooth transverse branches, tangent to the principal
directions and the asymptotic directions at \pp, if \pp\ is a
hyperbolic point of the surface;
\newline\noindent
$\bullet$ two smooth transverse branches, tangent to the principal
directions at \pp, if \pp\ is a non-umbilic elliptic point;
\newline\noindent
$\bullet$ three smooth transverse branches, if \pp\ is a generic
umbilic point;
\newline\noindent
$\bullet$ three branches tangent to the zero-curvature principal
direction, one of which is smooth and the other two ones have an
ordinary cusp, if \pp\ is an ordinary parabolic point;
\newline\noindent
$\bullet$ two smooth transverse branches, one of which is tangent to
the parabolic curve, if \pp\ is an elliptic cusp of Gauss;
\newline\noindent
$\bullet$ six smooth branches, five of which are tangent to the
parabolic curve and the other one is transverse, if \pp\ is a
hyperbolic cusp of Gauss.
\end{theorem}

Throughout this paper, we suppose that \pp\ is a generic umbilic of
$M$, that is, the quadratic part $q_0$ of $f_0$ is proportional to
$x^2+y^2$.
The genericity condition is that $q_0$ does not divide the cubic
part of $f_0$.
We may assume without loss of generality that the cubic part of $f_0$
equals $ax^3+bx^2y+axy^2+cy^3$ (see Lemma 1 of \cite{dggk}).
The genericity condition is then $b\not=c$.

Let us consider an $n$-parameter deformation
$$z=f(x,y;\tau)=f_0(x,y)+R(x,y;\tau)$$ of our surface $z=f_0(x,y)$,
where $R:\R^2\times\R^n\longrightarrow \R$ is a smooth mapping identically
vanishing at $\tau=0$.
Up to affine coordinate changes depending on
the deformation parameters, we may assume that the origin \pp\ is not moved
by the deformation and that at this point the tangent plane to the
surface is always $z=0$.
Thus, we assume that the linear part of $R$
vanishes for all $\tau$.

\begin{definition}
The {\it bifurcation diagram} of the vertex sets of the family of
surfaces $z=f$ is the germ at the origin of the closure of the set in
the $(n+1)$-space $\R^n\times\R=\{\tau,k\}$
formed by the elements $(\tau,k)$ such
that the level curve $f(\, \cdot\,;\tau)=k$ has (at least) a
degenerate vertex.

The {\it discriminant} of the vertex sets is the projection of the
singular locus of the  bifurcation diagram to the $n$-parameter space
$\{\tau\}$ by the mapping ``forgetting $k$''.
\end{definition}

The forms $x^2+y^2$, $x^2-y^2$ and $2xy$ form a basis of the
$3$-dimensional space $Q$ of the quadratic forms on
$\R^2=\{x,y\}$.
Therefore, the umbilic forms (i.e., proportional to $x^2+y^2$) span a
codimension $2$ subspace of this space $Q$
(according to the fact that a generic surface has only isolated
umbilics).

Let us denote by $T$ the affine plane $q_0+\langle
x^2-y^2,2xy\rangle_\R\subset Q$, which is transversal at $q_0$ to the
line $\langle x^2+y^2\rangle_\R$ of the umbilic forms in the space of
the quadratic forms $Q$.

Consider the mapping $F$ which associates to a parameter $\tau$ the
natural projection on $T\approx\R^2=\{\lambda,\mu\}$ of the
corresponding quadratic part of $f(\,\cdot\,;\tau)$.

\begin{definition}
The {\it rank} of the deformation $f$ of $f_0$ is the rank of the
derivative of $F:\R^n\longrightarrow\R^2$, at the parameters' origin.
\end{definition}

\begin{remark}
The small deformations of rank $0$ perturb the surface $z=f_0$ among
the surfaces having a generic umbilic at the origin.
Therefore they do not change the vertex set of the surface (up to
diffeomorphisms).
\end{remark}

Therefore, we will discuss only deformations of rank $1$ and
$2$.
Notice that a deformation has generically the maximum rank
possible (i.e., $1$ for $1$-parameter deformations and $2$ for
$(n\geq2)$-parameter deformations).

\begin{theorem}\label{thm:2}
For every $n\geq 2$, the germ at the origin of the discriminant of the
vertex set of any $n$-parameter rank $2$ deformation $z=f$ of the
surface $z=f_0$ is diffeomorphic to the germ at the origin of an
$(n-2)$-cylinder over the union of three transverse smooth curves on
the plane $T$, tangent to the lines $\mu=0,\pm\sqrt{3}\lambda$.
\end{theorem}

 The vertex set at a generic umbilic
point of a surface is the union of three transverse smooth branches
$C_1$, $C_2$ and $C_3$ (Theorem \ref{thm:0}). Each of these branches
can be seen as the union $C_i=C_i^+\cup C_i^-$ of two half-branches
issuing from the origin.

\begin{theorem}\label{thm:3}
Let $z=f$ be an $n$-parameter rank $2$ deformation of $z=f_0$.
Fix a ball centered at the origin $(x=0,y=0)$ of radius arbitrarily
small.
For any $\tau$ small enough outside the deformation's discriminant,
the vertex set of the surface $z=f(\, \cdot\, ;\tau)$ in this ball is
the union of a smooth branch not passing through the
origin and two smooth branches passing through the origin.
The first of these branches is a smooth $C^0$-small deformation of two
consecutive half-branches, say $C_1^+\cup C_2^+$, of the unperturbed
vertex set; it is disjoint to the other two, which are
$C^0$-small deformations of the pairs of non consecutive remaining
half-branches $C_1^-\cup C_3^-$ and $C_2^-\cup C_3^+$
(see Figure \ref{fig:vertex:0}).
\end{theorem}
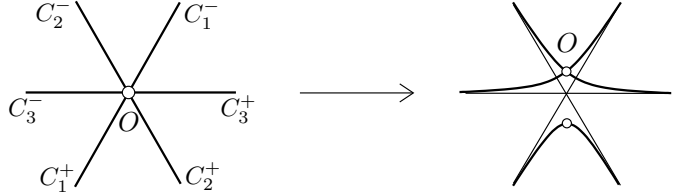
\begin{figure}[h]
 \centering
\scalebox{.95}{\input{fig-vertex-0.pstex_t}}
\caption{The vertex set at a generic umbilic point and its evolution
  under a generic small perturbation of the surface.}
 \label{fig:vertex:0}
\end{figure}

\begin{figure}[h]
\centering   \scalebox{.46}{\input{fig-vertex-4.pstex_t}}
 \caption{Vertex set discriminant of $2$-parameter rank $2$
   deformations of surfaces at generic umbilic points (middle picture)
   together with all corresponding perestroikas.}
 \label{fig:vertex:4}
\end{figure}
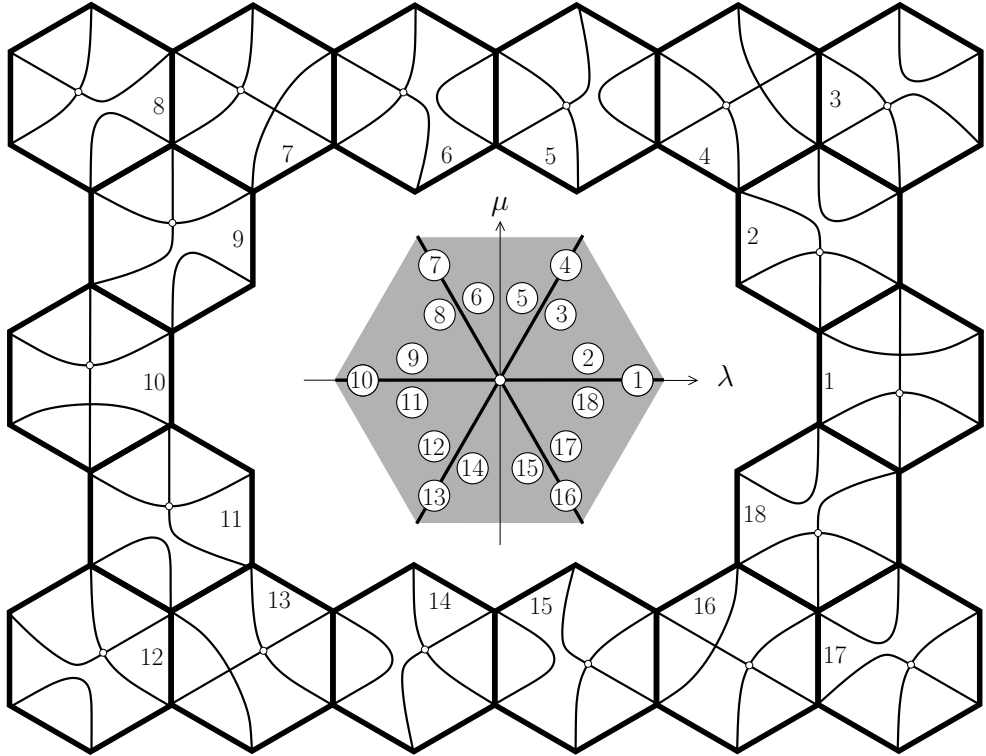

The projection on $T$ of the discriminant of any $2$-parameter rank $2$
deformation of a surface near a generic umbilic point is shown in
Figure \ref{fig:vertex:4}, together with the corresponding vertex sets.

\begin{remark}
Since generic umbilics are stable, the perturbed surface has a generic
umbilic close to that of the unperturbed surface.
Whenever the coordinates of the projection of the deformation $R$ on
the plane $T$ are not both vanishing, the umbilic slightly moves away
from the origin.
For small perturbations, the tangent plane to the perturbed surface at
the umbilic point is transversal to the planes $z=k$, so we cannot see
the umbilic in the diagram above as intersection of vertex set branches.
\end{remark}

\begin{example}
The vertex sets and some level lines of the surface $z=f$,
with function
$$f(x,y;\lambda,\mu)=x^2+y^2+x^3-y^3+\lambda(x^2-y^2)+2\mu xy \ , $$
are drawn in Figure \ref{fig:vertex:8} for the parameter values
$(\lambda,\mu)=(0,0), (1/10,0)$ and $(0,1/10)$.
\begin{figure}[h]
\centering
\includegraphics[width=0.22\textwidth]{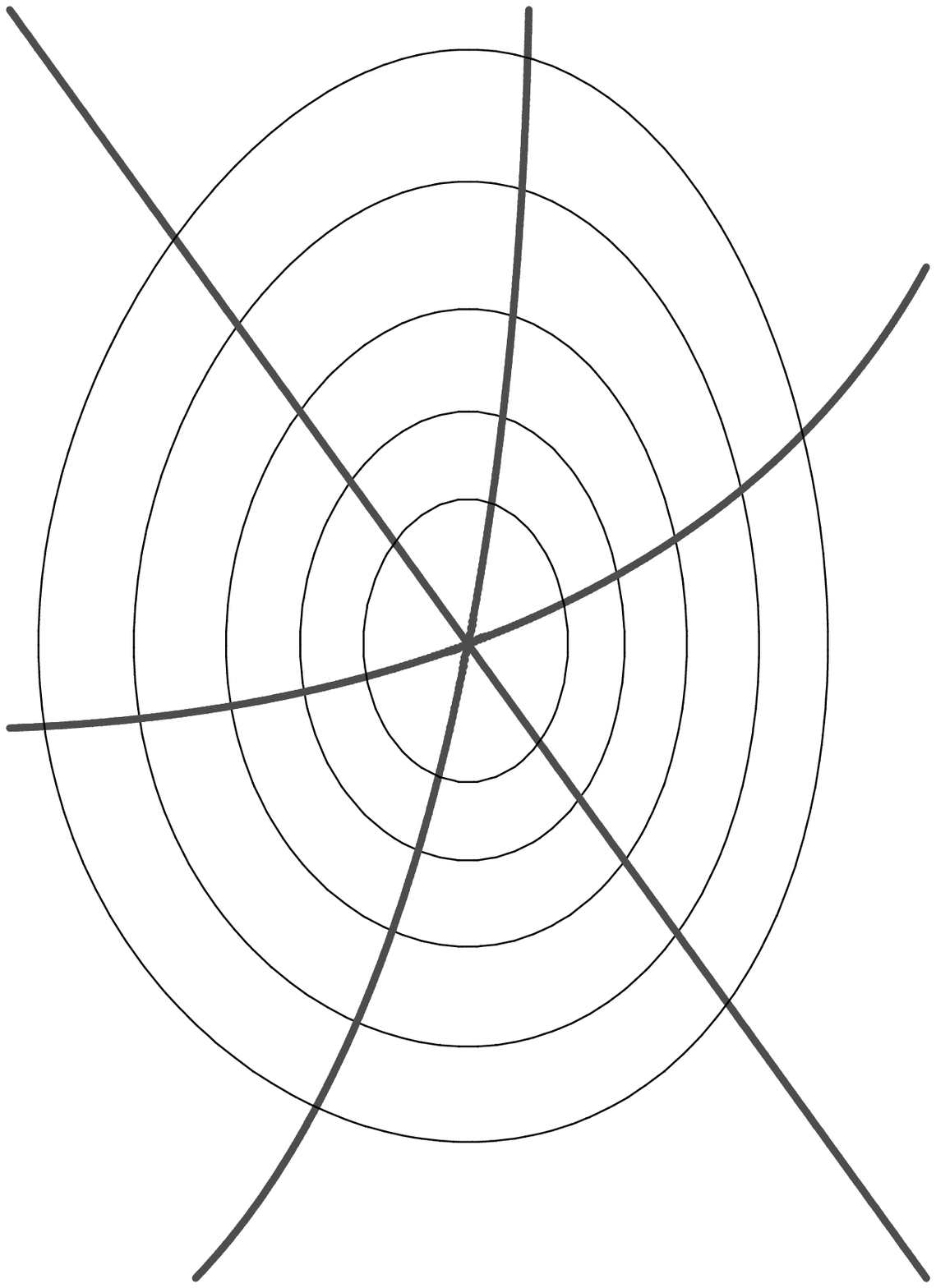} \qquad
\includegraphics[width=0.22\textwidth]{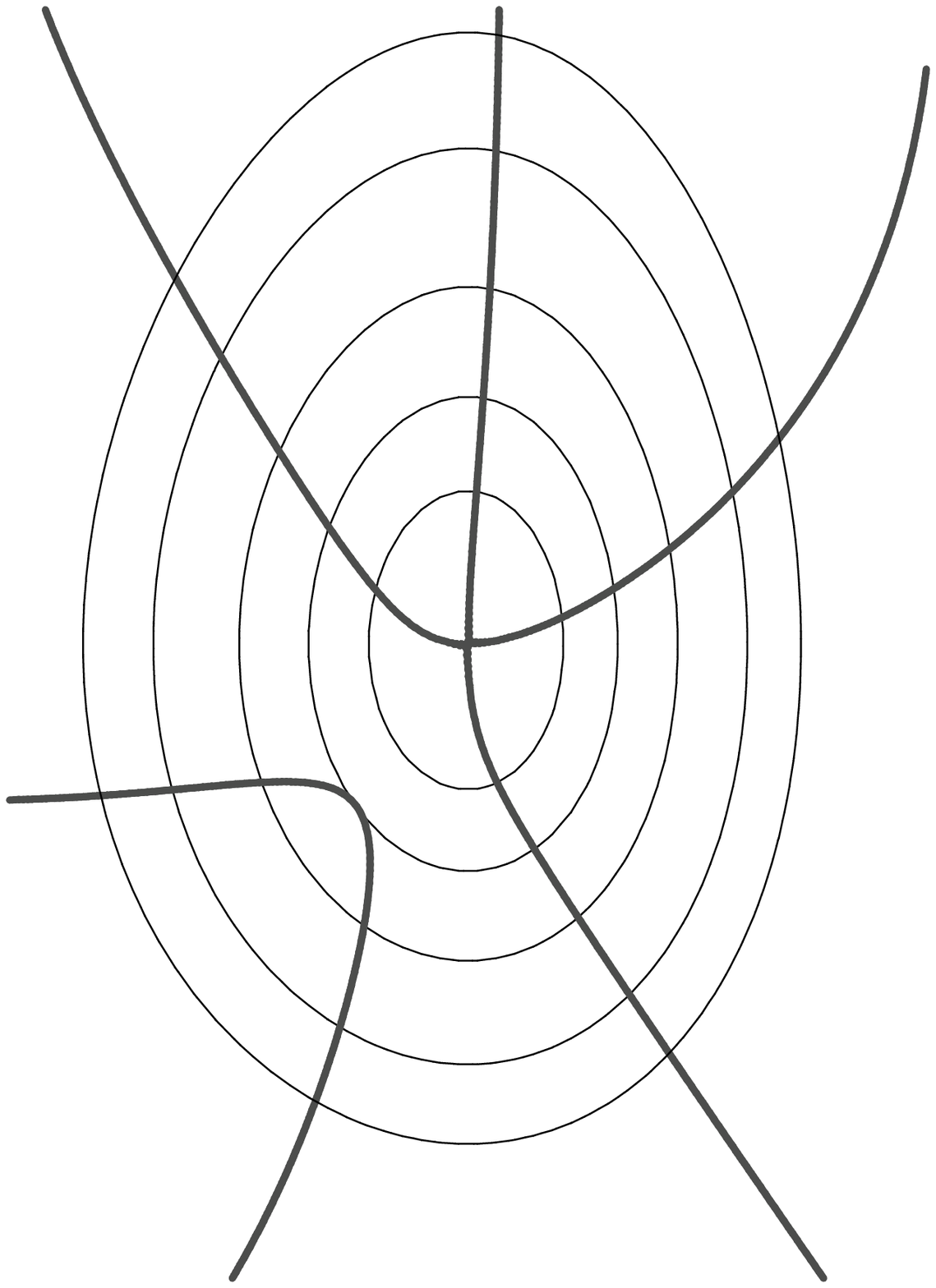} \qquad
\includegraphics[width=0.22\textwidth]{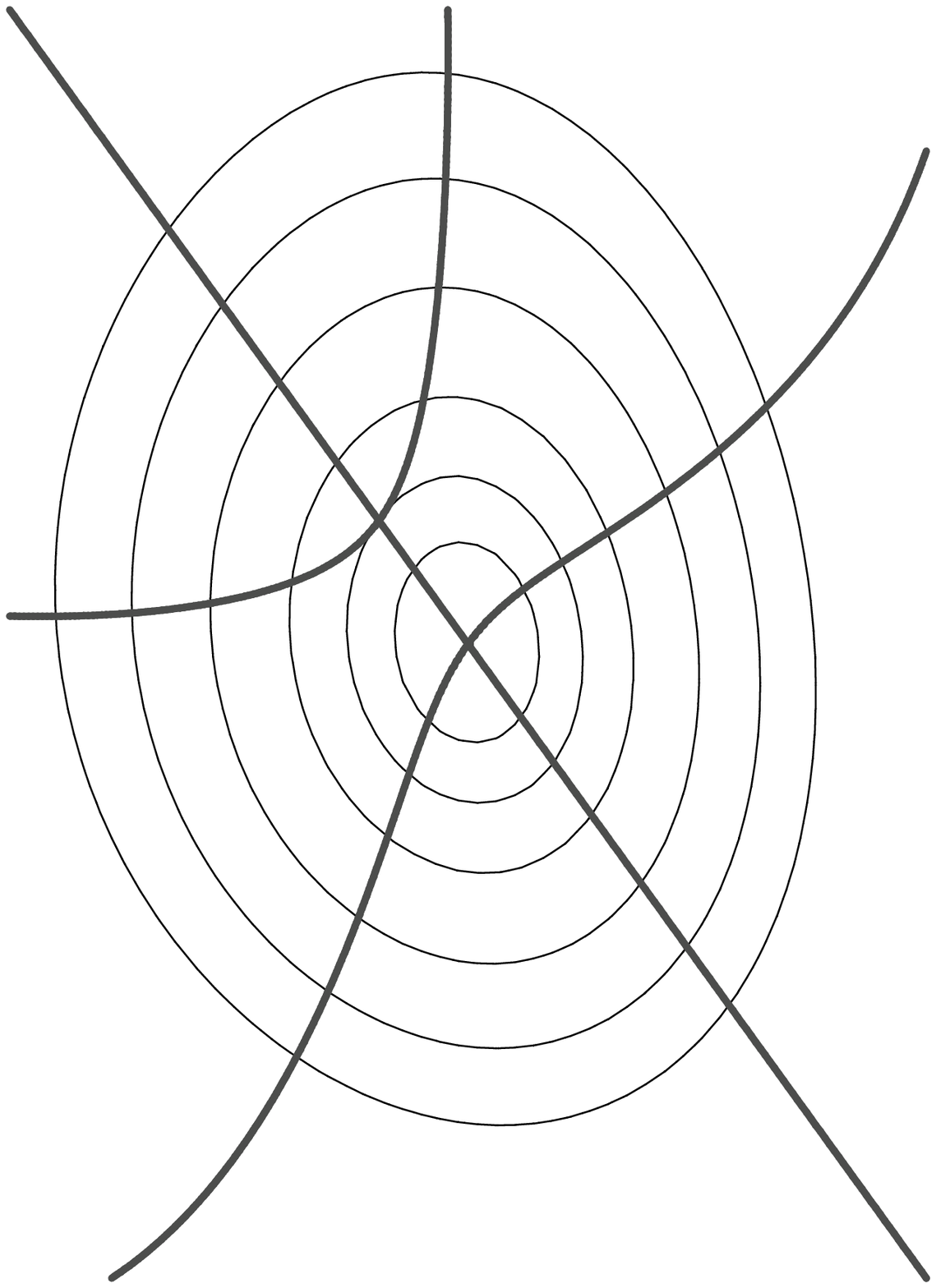}
\caption{Vertex sets of the surface $z=f$ for different values of the
  deforming parameters.}
\label{fig:vertex:8}
\end{figure}\end{example}

Proving a conjecture by V. Arnold,
R. Uribe-Vargas has shown that the bifurcation diagram of the vertex
sets of a generic $2$-parameter deformation of a surface at a generic
umbilic point is diffeomorphic to the {\em hypocycloidal cup} in the
$3$-space $\C\times\R$, parametrized by the mapping
$$(\phi,k)\mapsto \left( -2\sqrt{k}\,
  (5e^{-i\phi}+e^{5i\phi}),k\right)\ . $$

The curves above the hypocycloidal cup have $6$ vertices while the curves
below it have $4$ vertices.
The curves on the regular part of the hypocycloidal cup have a
$1$-degenerate vertex, while those on the semicubic cuspidal edge have
a $2$-degenerate vertex.

The proof, similar to that of the analoguous result on the nearby
problem on vanishing flattenings of spatial curves (see \cite{uribe}),
is based on Arnold's Lagrangian Collapse and Sturm-Hurwitz
Theorem (see \cite{arnoldlc}, \cite{ot}).
The relation between these two problems is explained in
Uribe-Vargas' comment to problem 1993-3 of \cite{arnoldp}.

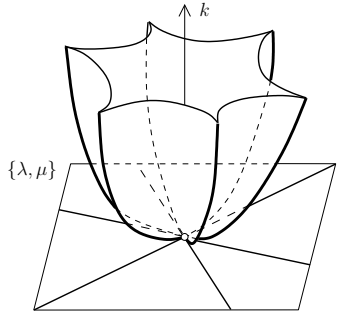
\begin{figure}[h]
  \centering
  \scalebox{.6}{\input{fig-vertex-6.pstex_t}}
  \caption{Arnold--Uribe-Vargas' hypocycloidal cup.}
  \label{fig:vertex:6}
\end{figure}

The discriminant of a rank $2$ perturbation is hence the projection on
the parameter plane of the semicubical cuspidal edges of the
hypocycloidal cup (see figure \ref{fig:vertex:6}).

\begin{remark}
Consider an $n$-parameter rank $2$ deformation of $z=f_0$.
If $\tau$ is small enough and belongs to the discriminant's
complement, then there exists a level curve $f(\,\cdot\, ;\tau)=k^*$
having a $1$-degenerate vertex.
This level separates the level curves having $4$ vertices
from those having $6$ vertices.
The value $k^*$ is arbitrarily small, provided that $\tau$ is small
enough: it is actually of the order of $\lambda^2+\mu^2$, as follows
from the parameterisation of the hypocycloidal cup (where $\lambda$ and
$\mu$ are the coordinates of the projection of the deformation $R$ on
the plane $T$).
\end{remark}

Since rank $1$ deformations are obviously induced from rank $2$
deformations, they are completely described by the above theorems.
In the most interesting case of generic $1$-parameter deformations we
obtain immediatly the following result.

\begin{theorem}\label{thm:1}
Let us consider a deformation of a surface $z=f_0$,
depending on one parameter $t\in\R$.
Assume that the projection on $T$ of the deformation is a curve
transversal at $\lambda=\mu=0$ to the lines
$\mu=0,\pm\sqrt{3}\lambda$.
Label the branches in such a way that for $t<0$ the vertex set
origin-avoiding branch is a smooth $C^0$-small deformation of the two
consecutive half-branches $C_1^+\cup C_2^+$ of the unperturbed vertex set.
Then for $t>0$ the vertex set origin-avoiding branch
is a smooth $C^0$-small deformation of the two opposite consecutive
half-branches $C_1^-\cup C_2^-$.
\end{theorem}

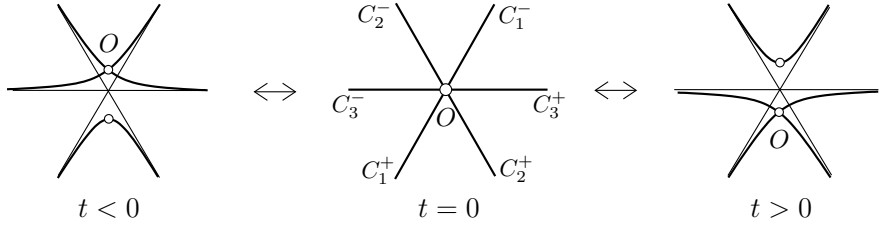
\begin{figure}[h]
\centering
\scalebox{.9}{\input{fig-vertex-7.pstex_t}}
\caption{Perestroika of the vertex set under a generic $1$-parameter
deformation of the surface at a generic umbilic.}
\label{fig:vertex:7}
\end{figure}

The perestroika of the vertex set under a generic $1$-parameter
deformation of the surface $z=f_0$ is illustrated in Figure \ref{fig:vertex:7}.

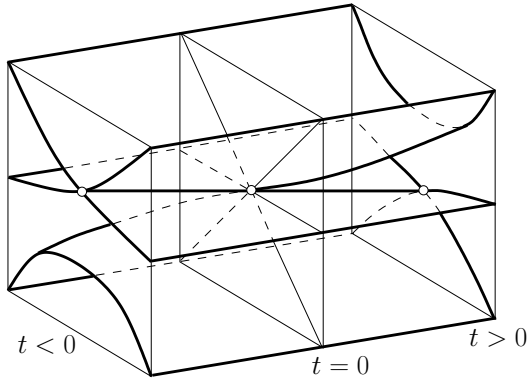
\begin{figure}[h]
\centering
\scalebox{.6}{\input{fig-vertex-5.pstex_t}}
\caption{Surface spanned by the deforming vertex sets on a generic
  $1$-parameter deformation. The slices of this surface by planes
  $t=const$ give the corresponding perestroika of the vertex set
  (compare Theorem \ref{thm:1} and Fig. \ref{fig:vertex:7}.)}
\label{fig:vertex:5}
\end{figure}

The germ at the origin of the surfaces formed in the $3$-space
$\{x,y,t\}$ by the vertex sets along generic $1$-parameter
deformations of the surface is depicted in figure \ref{fig:vertex:5}.

%*****************************************
\section{Proofs and further discussions}\label{chap:proofs}
%*****************************************

Let us consider as above a smooth $n$-parameter family of functions
$$f=f_0+R:\R^2\times\R^n\longrightarrow \R \ ,$$
with quadratic part equal to $x^2+y^2$ and cubic part equal to
$ax^3+bx^2y+axy^2+cy^3,$ where $b\neq c$ by the genericity
assumption.

To any fixed value $\tau$ of
the deforming parameters corresponds the vertex set of the perturbed
surface at the origin.
This vertex set is the zero level of a smooth function
$V_f(\, \cdot\,;\tau)$, depending smoothly on the parameter value
$\tau$.
Therefore we obtain a smooth $n$-parameter family of
{\it vertex set functions} $V_f:\R^2\times\R^n\longrightarrow\R$.
These functions are computed from $f$ by the following explicit formula,
given in \cite{dg}:
\begin{align}\label{eq:V}
V_f=&(f_x^2+f_y^2)(f_x^3f_{yyy}-3f_x^2 f_yf_{xyy}+3 f_x f_y^2
f_{xxy}-f_y^3f_{xxx})+ \nonumber\\
&\ \  +3f_x f_y \big(f_y^2 f_{xx}^2-f_x^2f_{yy}^2+ (f_x^2-f_y^2)
(f_{xx}f_{yy}+2f_{xy}^2) \big)+\nonumber \\
&\ \ +3 f_{xy} \big(f_{xx} f_y^4-3f_x^2f_y^2(f_{xx}-f_{yy})-
f_{yy}f_x^4\big) \ .
\end{align}

In order to describe geometrically these vertex sets, we consider the
natural equivalence relation acting on these functions.
Since we are interested on the zero levels of these functions, the
relevant equivalence relation is a version of the usual
$V$-equivalence, preserving the distinguished role of the deformation
parameters.

Consider the natural structure of trivial fiber bundle on
$\R^{2}\times\R^n$, defined by the natural projection
$\pi(x,y;\tau)=\tau$ onto the parameter space.

\begin{definition}
Two functions $F,G:\R^{2}\times\R^n\longrightarrow\R$ are
{\em $V^*$-equivalent} (or fibered $V$-equivalent) if there exist
diffeomorphisms
$$\Phi=(\Phi_1,\Phi_2):\R^2\times\R^n\longrightarrow\R^2\times\R^n \ ,
\quad \psi:\R\longrightarrow\R \ , $$
with $\psi(0)=0$,
making commutative the following diagram:
$$\begin{CD}\R^n @<\pi<< \R^2\times\R^n@>F>> \R  \\
@V\Phi_2 VV  @V\Phi VV     @VV\psi V         \\
\R^n @<\pi<< \R^2\times\R^n@>G>> \R \\
\end{CD}$$
\end{definition}

Notice that the left part of the diagram just means that $\Phi$ is a
fibered diffeomorphism of the total space of the fiber bundle
$\R^2\times\R^{n}\longrightarrow\R^n$.
In particular, each $\Phi_\tau(\cdot):=\Phi_1(\cdot;\tau)$ is a
diffeomorphism of each fiber $\R^2$.

A similar definition of $V^*$-equivalence holds for germs.

\begin{proposition}\label{prop:0}
Suppose that $F, G:\R^{2}\times\R^n\longrightarrow\R$ are $V^*$-equivalent.
Then the zero level sets of the functions $F(\cdot;\tau)$ and
$G\big(\cdot,\Phi_2(\tau)\big)$ are diffeomorphic.
\end{proposition}

\begin{proof}
Let $X_\tau$ be the zero level set of $F(\cdot;\tau)$.
Hence
$$X_\tau=\{(x,y):F(x,y;\tau)=0\}=\{(x,y):\psi^{-1}\circ G\circ \Phi
(x,y;\tau)=0\} \ . $$
Since $\psi^{-1}(t)=0$ if and only if $t=0$, we have
$$X_\tau=\{(x,y):G\big(\Phi_\tau(x,y);\Phi_2(\tau)\big)=0\} \ . $$
Therefore the zero level set of $G(\cdot;\Phi_2(\tau))$ is
$\Phi_\tau(X_\tau)$.
\end{proof}

$V^*$-equivalence and Proposition \ref{prop:0} allow us to replace
the qualitative study of the vertex set by the corresponding study
of a diffeomorphic curve $\tilde V_f=0$ (which is not necessarily
the vertex set of some surface).

Let us consider an $n$-parameter rank $2$ deformation $z=f_0+R$ of the
surface $z=f_0$.
Then $n\geq 2$ and, by the implicit function theorem,
the perturbing term is (up to a coordinate change
in the parameter space) of the form
$$R(\tau)= \lambda\, (x^2-y^2)+2\mu\, xy +R_3(\tau)\ , $$
where $\lambda$ and $\mu$ are two distinguished parameters among the
$\tau$'s and, as a function of $x,y$, $R_3(\tau)$ has a $2$-jet which
is identically zero.

\begin{lemma}\label{lemma:1}
The vertex set function $V_f$ is $V^*$-equivalent to a function
$\tilde V_f$ such that
$$\tilde V_4 = (1-\lambda^2-\mu^2)^2\,\big(x^2+y^2+\lambda\,
(x^2-y^2)+2\mu\, xy\big)\, \big(2\lambda\, xy+\mu\, (y^2-x^2)\big) \ , $$
$$\tilde V_5\mid_{\tau=0}=x(x^2+y^2)(x^2-3y^2)\ , \quad \tilde
V_6\mid_{\tau=0}=0 \ , $$
where $\tilde V_f=\tilde V_4+\tilde V_5+\tilde V_6+\dots$ is the
homogeneous expansion of $\tilde V$ (each $\tilde V_i$ is a
homogeneous polynomial of degree $i$ in the variables $x$ and $y$,
whose coefficients depend on the parameters $\tau$).
\end{lemma}

\begin{proof}
Compute $V_f$ in terms of $f$ by the explicit
expression \eqref{eq:V}, and consider the $V^*$-equivalent function
$$\tilde V_f:=\frac{(c-b)^4}{192}\, V_f +q_1\, x V_f + q_2\ y V_f \ , $$
for some coefficients $q_i$ (depending on the parameters).
Consider a coordinate change of the form
$$\left(
\begin{matrix} x\\y\end{matrix}\right)\mapsto
\left(\begin{matrix}
\frac{1}{c-b}\, x +p_0\, x^2+p_1\, xy+p_2\, y^2 \\
\frac{1}{c-b}\, y +p_3\, x^2 \\
\end{matrix}\right) \ ,
$$
where the coefficients $p_i$ depend smoothly on the parameters.
One directly check that for a suitable choice of these coefficients
$p$, $q$ (depending also on the quartic part of $f$), this function
fullfils the required conditions.
Notice that $b\not=c$, since the umbilic is
generic.
\end{proof}

\begin{remark}
$\tilde V_4$ depends only on the parameters $\lambda$ and
$\mu$.
In particular, the coefficients of $x^4$ and $y^4$ of $\tilde V_4$
are respectively
$$-\mu\, (1+\lambda)\, (1-\lambda^2-\mu^2)^2 \quad \text{and} \quad
\mu\, (1-\lambda)\, (1-\lambda^2-\mu^2)^2  \ ,$$
so near $\tau=0$ they both vanish if and only if $\mu=0$.
\end{remark}

\begin{lemma}\label{lemma:2}
For every small enough value of the parameter deformation $\tau$, such
that $(\lambda,\mu)\not=(0,0)$, the vertex set of the perturbed
surface $z=f$ has two real smooth branches passing through the
origin.
They are tangent to the lines
$$\mu\, y=(-\lambda\pm \sqrt{\lambda^2+\mu^2})\, x  \, $$
for $\mu\not=0$ and to the lines $x=0$, $y=0$ for $\mu=0$ and
$\lambda\not=0$. For $\lambda=\mu=0$, the vertex set has three real
smooth branches passing through the origin, tangent to the lines
$x=0$ and $x=\pm\sqrt{3}\, y$.
\end{lemma}

\begin{proof}
We shall use the Newton diagrams of the vertex set functions for the
different fixed values of the parameters.
The relevant diagrams are shown in figure \ref{fig:vertex:1}.
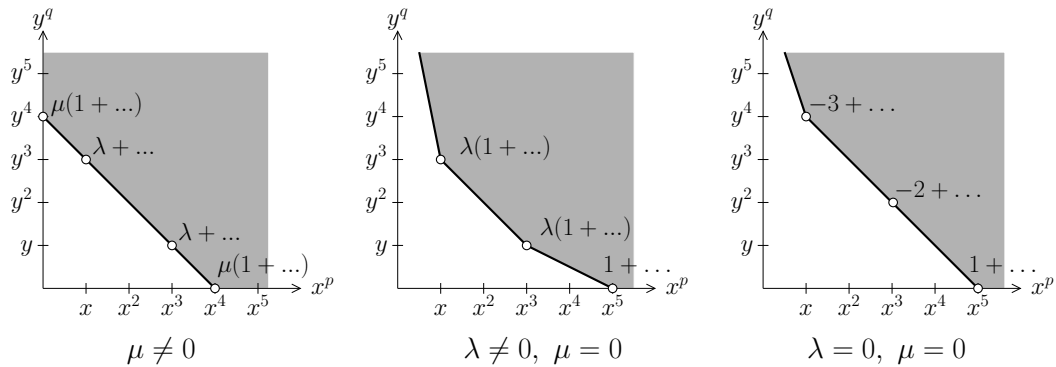
\begin{figure}[h]
  \centering
  \scalebox{.45}{\input{fig-vertex-1.pstex_t}}
  \caption{Newton diagrams of $\tilde V_f$.}
  \label{fig:vertex:1}
\end{figure}

Suppose first $\mu\not=0$.
Then $\tilde V_f$ is quasi-homogeneous, with principal part $\tilde
V_4$ (given in Lemma \ref{lemma:2}).
Solving $\tilde V_4=0$ we get:
$$y=\dfrac{-\lambda\pm \sqrt{\lambda^2+\mu^2}}{\mu}\, x \, \quad
y=\dfrac{\mu\pm \sqrt{\lambda^2+\mu^2-1}}{\lambda-1}\, x \ .$$
Hence, for $\lambda$ and $\mu$ small enough, $\tilde V_f=0$ has two
real smooth branches  passing through the origin, whose tangent
directions are given by the first two lines above.
The vertex set has also two smooth complex conjugate branches (tangent
to the second two lines above), whose
only real point is the origin (they appear also in the forthcoming
cases, but we do not insist about that).

Suppose now $\mu=0$, $\lambda\not=0$.
In this case the vertex set has two real smooth branches with equation
$$2\lambda\big(1+o(\lambda)\big)\, xy+hot(x,y)=0\ , $$
which are therefore tangent at the origin to the
lines $x=0$ and $y=0$.

Finally, assume $\lambda=\mu=0$.
Then the vertex set has three smooth real branches, according to the
fact that under such a deformation the origin is still an umbilic of
our surface.
These branches are provided by an equation of the form
$$x(x^2+y^2)(x^2-3y^2)+hot(x,y)=0 \ ,$$
so they are tangent to the lines $x=0,\pm\sqrt{3}\, y$.
\end{proof}

In order to describe the perestroikas of the vertex sets of the family
of surfaces $z=f(x,y)$, we first focalize on the case $\mu= 0$.

Let us fix a ball centered at the origin, of radius arbitrary
small.
For $\tau=0$, the real smooth branches of the vertex set
meet the ball's boundary at $6$ points, close to the vertex of the
regular hexagon inscribed in the ball and symmetric with respect to
$x=0$.
When $\tau$ varies, these six points slightly move along the
ball's boundary, provided that the variation is small
enough.

\begin{proposition}
Assume $\mu=0, \lambda\not=0$.
The branch of the vertex set, which is tangent to the line $y=0$, is
parabolic near the origin, namely it has a second order tangency with
a parabola of the form
$$2\lambda\, y = (-1+hot(\tau)) \,x^2 \ ,  \qquad \text{(for
  $\tau\rightarrow 0$)} \ . $$
The second vertex set branch has a second order tangency with its
tangent line $x=0$ at the origin.
\end{proposition}

\begin{proof}
The vertex set branch which is tangent to $y=0$ at the origin is of
the form $y=h(x)=C \,x^2+o(x^2)$, for some ``constant'' $C$
depending smoothly on the value of the parameter $\tau$.
Putting that into the expression of $\tilde V_f$, we get
(according to figure \ref{fig:vertex:1})
$$\tilde
V_f(x,h(x);\lambda,\mu=0)=
\big(1+2\lambda(1+hot(\tau))\ C+hot(\tau)\big)\, x^5+o(x^5)
=0\ .$$
Hence, $C=(-1+hot(\tau))/2\lambda$.

Similarly, we get the expression $x=o(y^2)$ for the other vertex set
smooth branch passing through the origin.
\end{proof}

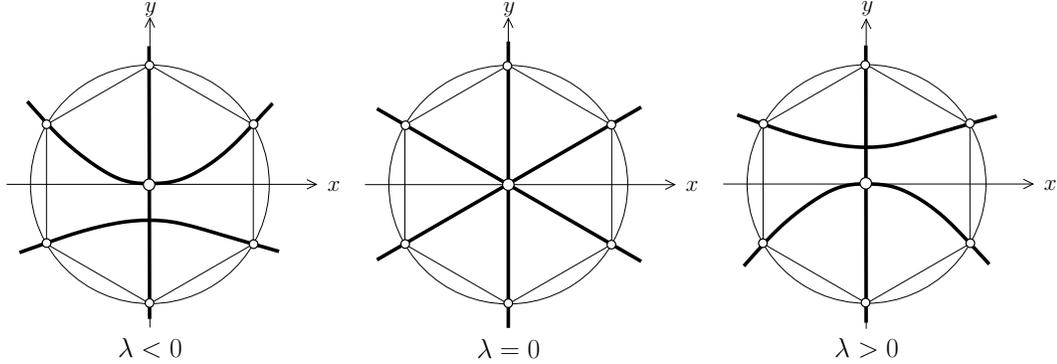
\begin{figure}[h]
\centering
\scalebox{.5} {\input{fig-vertex-2.pstex_t}}
\caption{Perestroika of the vertex set along the deformation $\mu=0$
  when $\lambda$ crosses $0$.}
\label{fig:vertex:2}
\end{figure}

This proposition allows us to describe the perestroika of the vertex
set along this deformation.
In particular, the convexity of this branch is toward the positive
direction  (resp. negative direction) of the $y$-axis whenever
$\lambda$ is negative (resp. positive).
The result is shown in figure
\ref{fig:vertex:2}.

Notice that the pattern of the third vertex set branch depicted in
figure \ref{fig:vertex:2} is the only possible (provided that the
deformation is small enough, as well as the radius of the ball inside
which we are looking the vertex set).
Indeed, for $\tau,k$ small enough, $k>0$, the level curve $f=k$ has at
most one $2$-degenerate vertex by Uribe-Vargas' theorem.

In particular, the vertex corresponding to the self-intersection of
the vertex set under our deformations is exactly $2$-degenerate.

The $6$-jet of the vertex set equation has an hexagonal symmetry which
simplifies our computations.
Indeed, one easily check the following fact.

\begin{proposition}
For every $n\in\Z$, the terms $\tilde V_4$, $\tilde V_5\mid_{\tau=0}$
and $\tilde V_6\mid_{\tau=0}$ are invariant under the simultaneous rotations
of angles $2n\pi/3$ and $-2n\pi/3$ in the $\{x,y\}$ plane and in the
$\{\lambda,\mu\}$ plane.
\end{proposition}

This remarkable hexagonal symmetry (related to the Sturm-Hurwitz
theorem, cf. \cite{uribe}) of the vertex set implies that the above
perestroika, occurring on the hypersurface  $\mu=0$ when $\lambda$
crosses $0$, also occurs
along two other smooth hyperplanes, whose projections on
the parameter plane $\{\lambda,\mu\}$ are the
lines $\mu=\pm\sqrt{3}\, \lambda$.

In particular, these three hypersurfaces are contained in the
deformation's discriminant.

A similar study of the relative positions of the vertex set branches
leads to the description of the perestroika of the vertex set along
the deforming curves $\mu=\alpha\,\lambda$, where
$\alpha\not=0,\pm\sqrt{3}$.

To describe the behaviour of the vertex set branches (in a small ball
centered at the origin), we shall compute the convexity of the
branches passing through the origin.

\begin{proposition}
Suppose $\tau$ small enough and such that $\mu=\alpha\,\lambda$,
$\alpha\not=0,\pm\sqrt{3}$.
Then the two smooth real branches of the vertex set passing through
the origin have there a second order tangency with two parabol\ae \ of
the form
\begin{equation}
\label{eq:2}
y=h(x)=\left(\dfrac{-1\pm\sqrt{1+\alpha^2}}{\alpha}+hot(\tau)\right) \,
x +A_\pm\, x^2+o(x^2) \ ,
\end{equation}
where the coefficients $A_\pm$, depending smoothly on $\tau$ and
$\alpha$, are
$$A_\pm=\dfrac{(\alpha^2+6)(\alpha^2+1)\mp 2\sqrt{1+\alpha^2}
  (2\alpha^2+3)}
{\lambda\alpha^2(1+\alpha^2)(-1\pm\sqrt{1+\alpha^2})}\ (1+hot(\tau)) \
  .  $$

In particular, $A_\pm\not=0$ for every $\lambda\not=0$ small enough.
\end{proposition}

\begin{proof}
By Lemma \ref{lemma:2}, the two real smooth branches of the vertex set
passing through the origin are of the form \eqref{eq:2}.
Replacing this expression in $\tilde V_f$ and using Lemma
\ref{lemma:1}, we get an explicit expression of the form
$C_\pm x^5+o(x^5)$.
The equation $C_\pm=0$ provides the above values of $A_\pm$.

The coefficient $A_-$ vanishes if
and only if $\alpha$ is a root of the polynomial
$$(\alpha^2+6)^2(\alpha^2+1)-4(2\alpha^2+3)^2= \alpha^4(\alpha^2-3)\ . $$
Thus $A_\pm\not=0$ for $\lambda$ small enough, under the hypothesis
of the proposition.
\end{proof}

\begin{figure}[h]
\centering
\scalebox{.5}{\input{fig-vertex-3.pstex_t}}
\caption{Perestroika of the vertex set along the deformation
  $\mu=\alpha\lambda$.}
\label{fig:vertex:3}
\end{figure}
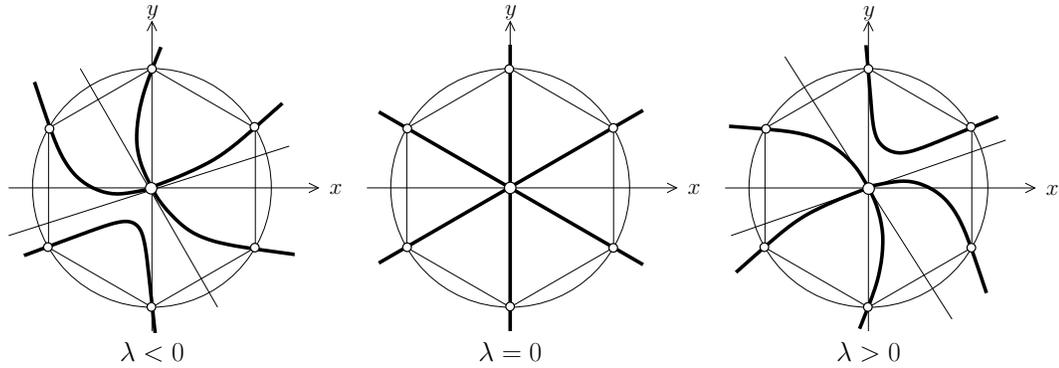

The resulting arrangement of the vertex set branches is shown in
figure \ref{fig:vertex:3} in the case $0<\alpha<\sqrt{3}$ (the other
cases can be deduced from this one via the hexagonal symmetry).
These arrangements are obtained using the convexity described in the above
proposition and the same argument used in the case $\mu=0$.
The key point here is that the tangent lines to the two branches
passing through the origin intersect the same sides of the symmetry
hexagon for every $0<\alpha<\sqrt{3}$.

In particular, the degenerate vertex on the origin-avoiding
branch is exactly $1$-degenerate.
These deformations belong to the complement of the
discriminant, for every $\tau$ small enough such that
$(\lambda,\mu)\not=(0,0)$.
Summing up the descriptions of these deformations, we have proved
Theorem \ref{thm:2} and Theorem \ref{thm:3}.

\section{Conclusion}\label{conclusion}
In this paper, we have discussed the geometry of the vertex sets at
isolated umbilic points of a generic surface evolving in an
$n$-parameter family of surfaces.  As a follow up, it would be very
interesting to extend this investigation to vertex sets at  higher
degenerate points (in the sense that the tangent plane has a higher
order contact with the surface) such as ordinary parabolic points or
cusps of Gauss, of $n$-parameter small deformations of generic
surfaces. For fixed (as opposed to continuous families of) generic
surfaces, these degenerations are a part of the work carried out in
\cite{dg}. Recall that the vertices we consider in this paper are
Euclidean invariants. A discussion of this subject from the
projective differential geometry point of view, is an interesting
question. A challenging problem would be, as brought to our
attention by V.I. Arnold, the extension of the work within this
paper to symplectic and contact geometries, namely the singularities
of the corresponding lagrangian map theory, where the vertices are
replaced by the singularities of the caustic enveloping the normals.
The theory in this context is still awaiting its creation. See
\cite{arnold96} for an extensive discussion.

\bigskip
\bigskip

\noindent{\sc Gianmarco Capitanio and Andr{\'e} Diatta}\\
Department of Mathematical Sciences,\\
Mathematics and Oceanography Building,\\
Peach Street,\\
Liverpool L69 7ZL \\
United Kingdom  \\
e-mail: capitani@math.jussieu.fr \quad(G. Capitanio)\\
e-mail: adiatta@liverpool.ac.uk \quad (A. Diatta)
\end{document}

%% file: fig-vertex-0.pstex_t
\begin{picture}(0,0)%
\includegraphics{fig-vertex-0.pstex}%
\end{picture}%
\setlength{\unitlength}{3947sp}%
\begingroup\makeatletter\ifx\SetFigFont\undefined%
\gdef\SetFigFont#1#2#3#4#5{%
  \reset@font\fontsize{#1}{#2pt}%
  \fontfamily{#3}\fontseries{#4}\fontshape{#5}%
  \selectfont}%
\fi\endgroup%
\begin{picture}(4562,1288)(609,-3421)
\put(4461,-2493){\makebox(0,0)[b]{\smash{{\SetFigFont{11}{13.2}{\rmdefault}{\mddefault}{\updefault}$O$}}}}
\put(886,-2914){\makebox(0,0)[b]{\smash{{\SetFigFont{10}{12.0}{\rmdefault}{\mddefault}{\updefault}$C_3^-$}}}}
\put(1092,-3372){\makebox(0,0)[b]{\smash{{\SetFigFont{10}{12.0}{\rmdefault}{\mddefault}{\updefault}$C_1^+$}}}}
\put(1564,-3015){\makebox(0,0)[b]{\smash{{\SetFigFont{11}{13.2}{\rmdefault}{\mddefault}{\updefault}$O$}}}}
\put(2055,-3360){\makebox(0,0)[b]{\smash{{\SetFigFont{10}{12.0}{\rmdefault}{\mddefault}{\updefault}$C_2^+$}}}}
\put(2288,-2915){\makebox(0,0)[b]{\smash{{\SetFigFont{10}{12.0}{\rmdefault}{\mddefault}{\updefault}$C_3^+$}}}}
\put(2048,-2291){\makebox(0,0)[b]{\smash{{\SetFigFont{10}{12.0}{\rmdefault}{\mddefault}{\updefault}$C_1^-$}}}}
\put(1066,-2285){\makebox(0,0)[b]{\smash{{\SetFigFont{10}{12.0}{\rmdefault}{\mddefault}{\updefault}$C_2^-$}}}}
\end{picture}%

%% file: fig-vertex-4.pstex_t
\begin{picture}(0,0)%
\includegraphics{fig-vertex-4.pstex}%
\end{picture}%
\setlength{\unitlength}{3947sp}%
\begingroup\makeatletter\ifx\SetFigFont\undefined%
\gdef\SetFigFont#1#2#3#4#5{%
  \reset@font\fontsize{#1}{#2pt}%
  \fontfamily{#3}\fontseries{#4}\fontshape{#5}%
  \selectfont}%
\fi\endgroup%
\begin{picture}(13396,10338)(-8257,-13788)
\put(-1493,-6323){\makebox(0,0)[b]{\smash{{\SetFigFont{25}{30.0}{\rmdefault}{\mddefault}{\updefault}$\mu$}}}}
\put(1582,-8723){\makebox(0,0)[b]{\smash{{\SetFigFont{25}{30.0}{\rmdefault}{\mddefault}{\updefault}$\lambda$}}}}
\put(392,-8758){\makebox(0,0)[b]{\smash{{\SetFigFont{20}{24.0}{\rmdefault}{\mddefault}{\updefault}1}}}}
\put(-288,-8448){\makebox(0,0)[b]{\smash{{\SetFigFont{20}{24.0}{\rmdefault}{\mddefault}{\updefault}2}}}}
\put(-648,-7848){\makebox(0,0)[b]{\smash{{\SetFigFont{20}{24.0}{\rmdefault}{\mddefault}{\updefault}3}}}}
\put(-605,-7175){\makebox(0,0)[b]{\smash{{\SetFigFont{20}{24.0}{\rmdefault}{\mddefault}{\updefault}4}}}}
\put(-1190,-7625){\makebox(0,0)[b]{\smash{{\SetFigFont{20}{24.0}{\rmdefault}{\mddefault}{\updefault}5}}}}
\put(-1813,-7619){\makebox(0,0)[b]{\smash{{\SetFigFont{20}{24.0}{\rmdefault}{\mddefault}{\updefault}6}}}}
\put(-2377,-7185){\makebox(0,0)[b]{\smash{{\SetFigFont{20}{24.0}{\rmdefault}{\mddefault}{\updefault}7}}}}
\put(-2316,-7858){\makebox(0,0)[b]{\smash{{\SetFigFont{20}{24.0}{\rmdefault}{\mddefault}{\updefault}8}}}}
\put(-2676,-8451){\makebox(0,0)[b]{\smash{{\SetFigFont{20}{24.0}{\rmdefault}{\mddefault}{\updefault}9}}}}
\put(-3385,-8751){\makebox(0,0)[b]{\smash{{\SetFigFont{20}{24.0}{\rmdefault}{\mddefault}{\updefault}10}}}}
\put(-2697,-9058){\makebox(0,0)[b]{\smash{{\SetFigFont{20}{24.0}{\rmdefault}{\mddefault}{\updefault}11}}}}
\put(-2407,-9665){\makebox(0,0)[b]{\smash{{\SetFigFont{20}{24.0}{\rmdefault}{\mddefault}{\updefault}12}}}}
\put(-2396,-10345){\makebox(0,0)[b]{\smash{{\SetFigFont{20}{24.0}{\rmdefault}{\mddefault}{\updefault}13}}}}
\put(-1872,-9958){\makebox(0,0)[b]{\smash{{\SetFigFont{20}{24.0}{\rmdefault}{\mddefault}{\updefault}14}}}}
\put(-1127,-9959){\makebox(0,0)[b]{\smash{{\SetFigFont{20}{24.0}{\rmdefault}{\mddefault}{\updefault}15}}}}
\put(-612,-10344){\makebox(0,0)[b]{\smash{{\SetFigFont{20}{24.0}{\rmdefault}{\mddefault}{\updefault}16}}}}
\put(-606,-9672){\makebox(0,0)[b]{\smash{{\SetFigFont{20}{24.0}{\rmdefault}{\mddefault}{\updefault}17}}}}
\put(-306,-9052){\makebox(0,0)[b]{\smash{{\SetFigFont{20}{24.0}{\rmdefault}{\mddefault}{\updefault}18}}}}
\put(3001,-8766){\makebox(0,0)[b]{\smash{{\SetFigFont{20}{24.0}{\rmdefault}{\mddefault}{\updefault}1}}}}
\put(1306,-5684){\makebox(0,0)[b]{\smash{{\SetFigFont{20}{24.0}{\rmdefault}{\mddefault}{\updefault}4}}}}
\put(-4394,-5664){\makebox(0,0)[b]{\smash{{\SetFigFont{20}{24.0}{\rmdefault}{\mddefault}{\updefault}7}}}}
\put(-4509,-11756){\makebox(0,0)[b]{\smash{{\SetFigFont{20}{24.0}{\rmdefault}{\mddefault}{\updefault}13}}}}
\put(-6215,-8781){\makebox(0,0)[b]{\smash{{\SetFigFont{20}{24.0}{\rmdefault}{\mddefault}{\updefault}10}}}}
\put(1959,-6782){\makebox(0,0)[b]{\smash{{\SetFigFont{20}{24.0}{\rmdefault}{\mddefault}{\updefault}2}}}}
\put(3085,-4908){\makebox(0,0)[b]{\smash{{\SetFigFont{20}{24.0}{\rmdefault}{\mddefault}{\updefault}3}}}}
\put(-804,-5691){\makebox(0,0)[b]{\smash{{\SetFigFont{20}{24.0}{\rmdefault}{\mddefault}{\updefault}5}}}}
\put(-2220,-5679){\makebox(0,0)[b]{\smash{{\SetFigFont{20}{24.0}{\rmdefault}{\mddefault}{\updefault}6}}}}
\put(-6153,-5004){\makebox(0,0)[b]{\smash{{\SetFigFont{20}{24.0}{\rmdefault}{\mddefault}{\updefault}8}}}}
\put(-5069,-6818){\makebox(0,0)[b]{\smash{{\SetFigFont{20}{24.0}{\rmdefault}{\mddefault}{\updefault}9}}}}
\put(-6250,-12539){\makebox(0,0)[b]{\smash{{\SetFigFont{20}{24.0}{\rmdefault}{\mddefault}{\updefault}12}}}}
\put(-2318,-11772){\makebox(0,0)[b]{\smash{{\SetFigFont{20}{24.0}{\rmdefault}{\mddefault}{\updefault}14}}}}
\put(-934,-11839){\makebox(0,0)[b]{\smash{{\SetFigFont{20}{24.0}{\rmdefault}{\mddefault}{\updefault}15}}}}
\put(1298,-11836){\makebox(0,0)[b]{\smash{{\SetFigFont{20}{24.0}{\rmdefault}{\mddefault}{\updefault}16}}}}
\put(3076,-12511){\makebox(0,0)[b]{\smash{{\SetFigFont{20}{24.0}{\rmdefault}{\mddefault}{\updefault}17}}}}
\put(1982,-10575){\makebox(0,0)[b]{\smash{{\SetFigFont{20}{24.0}{\rmdefault}{\mddefault}{\updefault}18}}}}
\put(-5164,-10651){\makebox(0,0)[b]{\smash{{\SetFigFont{20}{24.0}{\rmdefault}{\mddefault}{\updefault}11}}}}
\end{picture}%

%% file: fig-vertex-6.pstex_t
\begin{picture}(0,0)%
\includegraphics{fig-vertex-6.pstex}%
\end{picture}%
\setlength{\unitlength}{3947sp}%
\begingroup\makeatletter\ifx\SetFigFont\undefined%
\gdef\SetFigFont#1#2#3#4#5{%
  \reset@font\fontsize{#1}{#2pt}%
  \fontfamily{#3}\fontseries{#4}\fontshape{#5}%
  \selectfont}%
\fi\endgroup%
\begin{picture}(3195,3238)(6103,-5910)
\put(6104,-4478){\makebox(0,0)[b]{\smash{{\SetFigFont{12}{14.4}{\rmdefault}{\mddefault}{\updefault}$\{\lambda,\mu\}$}}}}
\put(7910,-2810){\makebox(0,0)[b]{\smash{{\SetFigFont{12}{14.4}{\rmdefault}{\mddefault}{\updefault}$k$}}}}
\end{picture}%

%% file: fig-vertex-7.pstex_t
\begin{picture}(0,0)%
\includegraphics{fig-vertex-7.pstex}%
\end{picture}%
\setlength{\unitlength}{3947sp}%
\begingroup\makeatletter\ifx\SetFigFont\undefined%
\gdef\SetFigFont#1#2#3#4#5{%
  \reset@font\fontsize{#1}{#2pt}%
  \fontfamily{#3}\fontseries{#4}\fontshape{#5}%
  \selectfont}%
\fi\endgroup%
\begin{picture}(6124,1570)(-1523,-3703)
\put(-789,-2493){\makebox(0,0)[b]{\smash{{\SetFigFont{11}{13.2}{\rmdefault}{\mddefault}{\updefault}$O$}}}}
\put(3890,-3167){\makebox(0,0)[b]{\smash{{\SetFigFont{11}{13.2}{\rmdefault}{\mddefault}{\updefault}$O$}}}}
\put(886,-2914){\makebox(0,0)[b]{\smash{{\SetFigFont{10}{12.0}{\rmdefault}{\mddefault}{\updefault}$C_3^-$}}}}
\put(1092,-3372){\makebox(0,0)[b]{\smash{{\SetFigFont{10}{12.0}{\rmdefault}{\mddefault}{\updefault}$C_1^+$}}}}
\put(1564,-3015){\makebox(0,0)[b]{\smash{{\SetFigFont{11}{13.2}{\rmdefault}{\mddefault}{\updefault}$O$}}}}
\put(2055,-3360){\makebox(0,0)[b]{\smash{{\SetFigFont{10}{12.0}{\rmdefault}{\mddefault}{\updefault}$C_2^+$}}}}
\put(2288,-2915){\makebox(0,0)[b]{\smash{{\SetFigFont{10}{12.0}{\rmdefault}{\mddefault}{\updefault}$C_3^+$}}}}
\put(2048,-2291){\makebox(0,0)[b]{\smash{{\SetFigFont{10}{12.0}{\rmdefault}{\mddefault}{\updefault}$C_1^-$}}}}
\put(1066,-2285){\makebox(0,0)[b]{\smash{{\SetFigFont{10}{12.0}{\rmdefault}{\mddefault}{\updefault}$C_2^-$}}}}
\put(1583,-3647){\makebox(0,0)[b]{\smash{{\SetFigFont{12}{14.4}{\rmdefault}{\mddefault}{\updefault}$t=0$}}}}
\put(3902,-3649){\makebox(0,0)[b]{\smash{{\SetFigFont{12}{14.4}{\rmdefault}{\mddefault}{\updefault}$t>0$}}}}
\put(-790,-3646){\makebox(0,0)[b]{\smash{{\SetFigFont{12}{14.4}{\rmdefault}{\mddefault}{\updefault}$t<0$}}}}
\end{picture}%

%% file: fig-vertex-5.pstex_t
\begin{picture}(0,0)%
\includegraphics{fig-vertex-5.pstex}%
\end{picture}%
\setlength{\unitlength}{3947sp}%
\begingroup\makeatletter\ifx\SetFigFont\undefined%
\gdef\SetFigFont#1#2#3#4#5{%
  \reset@font\fontsize{#1}{#2pt}%
  \fontfamily{#3}\fontseries{#4}\fontshape{#5}%
  \selectfont}%
\fi\endgroup%
\begin{picture}(5504,3978)(3835,-5083)
\put(4276,-4786){\makebox(0,0)[b]{\smash{{\SetFigFont{17}{20.4}{\rmdefault}{\mddefault}{\updefault}$t<0$}}}}
\put(7351,-5011){\makebox(0,0)[b]{\smash{{\SetFigFont{17}{20.4}{\rmdefault}{\mddefault}{\updefault}$t=0$}}}}
\put(9001,-4711){\makebox(0,0)[b]{\smash{{\SetFigFont{17}{20.4}{\rmdefault}{\mddefault}{\updefault}$t>0$}}}}
\end{picture}%

%% file: fig-vertex-1.pstex_t
\begin{picture}(0,0)%
\includegraphics{fig-vertex-1.pstex}%
\end{picture}%
\setlength{\unitlength}{3947sp}%
\begingroup\makeatletter\ifx\SetFigFont\undefined%
\gdef\SetFigFont#1#2#3#4#5{%
  \reset@font\fontsize{#1}{#2pt}%
  \fontfamily{#3}\fontseries{#4}\fontshape{#5}%
  \selectfont}%
\fi\endgroup%
\begin{picture}(14943,4965)(6104,-5935)
\put(17401,-5236){\makebox(0,0)[b]{\smash{{\SetFigFont{20}{24.0}{\rmdefault}{\mddefault}{\updefault}$x$}}}}
\put(18001,-5236){\makebox(0,0)[b]{\smash{{\SetFigFont{20}{24.0}{\rmdefault}{\mddefault}{\updefault}$x^2$}}}}
\put(18601,-5236){\makebox(0,0)[b]{\smash{{\SetFigFont{20}{24.0}{\rmdefault}{\mddefault}{\updefault}$x^3$}}}}
\put(19201,-5236){\makebox(0,0)[b]{\smash{{\SetFigFont{20}{24.0}{\rmdefault}{\mddefault}{\updefault}$x^4$}}}}
\put(19801,-5236){\makebox(0,0)[b]{\smash{{\SetFigFont{20}{24.0}{\rmdefault}{\mddefault}{\updefault}$x^5$}}}}
\put(7351,-5236){\makebox(0,0)[b]{\smash{{\SetFigFont{20}{24.0}{\rmdefault}{\mddefault}{\updefault}$x$}}}}
\put(7951,-5236){\makebox(0,0)[b]{\smash{{\SetFigFont{20}{24.0}{\rmdefault}{\mddefault}{\updefault}$x^2$}}}}
\put(8551,-5236){\makebox(0,0)[b]{\smash{{\SetFigFont{20}{24.0}{\rmdefault}{\mddefault}{\updefault}$x^3$}}}}
\put(9151,-5236){\makebox(0,0)[b]{\smash{{\SetFigFont{20}{24.0}{\rmdefault}{\mddefault}{\updefault}$x^4$}}}}
\put(9751,-5236){\makebox(0,0)[b]{\smash{{\SetFigFont{20}{24.0}{\rmdefault}{\mddefault}{\updefault}$x^5$}}}}
\put(12301,-5236){\makebox(0,0)[b]{\smash{{\SetFigFont{20}{24.0}{\rmdefault}{\mddefault}{\updefault}$x$}}}}
\put(12901,-5236){\makebox(0,0)[b]{\smash{{\SetFigFont{20}{24.0}{\rmdefault}{\mddefault}{\updefault}$x^2$}}}}
\put(13501,-5236){\makebox(0,0)[b]{\smash{{\SetFigFont{20}{24.0}{\rmdefault}{\mddefault}{\updefault}$x^3$}}}}
\put(14101,-5236){\makebox(0,0)[b]{\smash{{\SetFigFont{20}{24.0}{\rmdefault}{\mddefault}{\updefault}$x^4$}}}}
\put(14701,-5236){\makebox(0,0)[b]{\smash{{\SetFigFont{20}{24.0}{\rmdefault}{\mddefault}{\updefault}$x^5$}}}}
\put(16576,-4336){\makebox(0,0)[b]{\smash{{\SetFigFont{20}{24.0}{\rmdefault}{\mddefault}{\updefault}$y$}}}}
\put(16501,-1936){\makebox(0,0)[b]{\smash{{\SetFigFont{20}{24.0}{\rmdefault}{\mddefault}{\updefault}$y^5$}}}}
\put(16501,-2536){\makebox(0,0)[b]{\smash{{\SetFigFont{20}{24.0}{\rmdefault}{\mddefault}{\updefault}$y^4$}}}}
\put(16501,-3136){\makebox(0,0)[b]{\smash{{\SetFigFont{20}{24.0}{\rmdefault}{\mddefault}{\updefault}$y^3$}}}}
\put(16501,-3736){\makebox(0,0)[b]{\smash{{\SetFigFont{20}{24.0}{\rmdefault}{\mddefault}{\updefault}$y^2$}}}}
\put(16801,-1186){\makebox(0,0)[b]{\smash{{\SetFigFont{20}{24.0}{\rmdefault}{\mddefault}{\updefault}$y^q$}}}}
\put(20701,-4936){\makebox(0,0)[b]{\smash{{\SetFigFont{20}{24.0}{\rmdefault}{\mddefault}{\updefault}$x^p$}}}}
\put(18526,-5836){\makebox(0,0)[b]{\smash{{\SetFigFont{25}{30.0}{\rmdefault}{\mddefault}{\updefault}$\lambda=0,\ \mu=0$}}}}
\put(6526,-4336){\makebox(0,0)[b]{\smash{{\SetFigFont{20}{24.0}{\rmdefault}{\mddefault}{\updefault}$y$}}}}
\put(6451,-1936){\makebox(0,0)[b]{\smash{{\SetFigFont{20}{24.0}{\rmdefault}{\mddefault}{\updefault}$y^5$}}}}
\put(6451,-2536){\makebox(0,0)[b]{\smash{{\SetFigFont{20}{24.0}{\rmdefault}{\mddefault}{\updefault}$y^4$}}}}
\put(6451,-3136){\makebox(0,0)[b]{\smash{{\SetFigFont{20}{24.0}{\rmdefault}{\mddefault}{\updefault}$y^3$}}}}
\put(6451,-3736){\makebox(0,0)[b]{\smash{{\SetFigFont{20}{24.0}{\rmdefault}{\mddefault}{\updefault}$y^2$}}}}
\put(6751,-1186){\makebox(0,0)[b]{\smash{{\SetFigFont{20}{24.0}{\rmdefault}{\mddefault}{\updefault}$y^q$}}}}
\put(10651,-4936){\makebox(0,0)[b]{\smash{{\SetFigFont{20}{24.0}{\rmdefault}{\mddefault}{\updefault}$x^p$}}}}
\put(8401,-5836){\makebox(0,0)[b]{\smash{{\SetFigFont{25}{30.0}{\rmdefault}{\mddefault}{\updefault}$\mu\not=0$}}}}
\put(11476,-4336){\makebox(0,0)[b]{\smash{{\SetFigFont{20}{24.0}{\rmdefault}{\mddefault}{\updefault}$y$}}}}
\put(11401,-1936){\makebox(0,0)[b]{\smash{{\SetFigFont{20}{24.0}{\rmdefault}{\mddefault}{\updefault}$y^5$}}}}
\put(11401,-2536){\makebox(0,0)[b]{\smash{{\SetFigFont{20}{24.0}{\rmdefault}{\mddefault}{\updefault}$y^4$}}}}
\put(11401,-3136){\makebox(0,0)[b]{\smash{{\SetFigFont{20}{24.0}{\rmdefault}{\mddefault}{\updefault}$y^3$}}}}
\put(11401,-3736){\makebox(0,0)[b]{\smash{{\SetFigFont{20}{24.0}{\rmdefault}{\mddefault}{\updefault}$y^2$}}}}
\put(11701,-1186){\makebox(0,0)[b]{\smash{{\SetFigFont{20}{24.0}{\rmdefault}{\mddefault}{\updefault}$y^q$}}}}
\put(15601,-4936){\makebox(0,0)[b]{\smash{{\SetFigFont{20}{24.0}{\rmdefault}{\mddefault}{\updefault}$x^p$}}}}
\put(13726,-5836){\makebox(0,0)[b]{\smash{{\SetFigFont{25}{30.0}{\rmdefault}{\mddefault}{\updefault}$\lambda\not=0,\ \mu=0$}}}}
\put(19276,-3586){\makebox(0,0)[b]{\smash{{\SetFigFont{20}{24.0}{\rmdefault}{\mddefault}{\updefault}$-2+\dots$}}}}
\put(20176,-4636){\makebox(0,0)[b]{\smash{{\SetFigFont{20}{24.0}{\rmdefault}{\mddefault}{\updefault}$1+\dots$}}}}
\put(15076,-4636){\makebox(0,0)[b]{\smash{{\SetFigFont{20}{24.0}{\rmdefault}{\mddefault}{\updefault}$1+\dots$}}}}
\put(7876,-2986){\makebox(0,0)[b]{\smash{{\SetFigFont{20}{24.0}{\rmdefault}{\mddefault}{\updefault}$\lambda+...$}}}}
\put(9826,-4636){\makebox(0,0)[b]{\smash{{\SetFigFont{20}{24.0}{\rmdefault}{\mddefault}{\updefault}$\mu(1+...)$}}}}
\put(7501,-2386){\makebox(0,0)[b]{\smash{{\SetFigFont{20}{24.0}{\rmdefault}{\mddefault}{\updefault}$\mu(1+...)$}}}}
\put(9076,-4186){\makebox(0,0)[b]{\smash{{\SetFigFont{20}{24.0}{\rmdefault}{\mddefault}{\updefault}$\lambda+...$}}}}
\put(14326,-4111){\makebox(0,0)[b]{\smash{{\SetFigFont{20}{24.0}{\rmdefault}{\mddefault}{\updefault}$\lambda(1+...)$}}}}
\put(13201,-2986){\makebox(0,0)[b]{\smash{{\SetFigFont{20}{24.0}{\rmdefault}{\mddefault}{\updefault}$\lambda(1+...)$}}}}
\put(18076,-2386){\makebox(0,0)[b]{\smash{{\SetFigFont{20}{24.0}{\rmdefault}{\mddefault}{\updefault}$-3+\dots$}}}}
\end{picture}%

%% file: fig-vertex-2.pstex_t
\begin{picture}(0,0)%
\includegraphics{fig-vertex-2.pstex}%
\end{picture}%
\setlength{\unitlength}{3947sp}%
\begingroup\makeatletter\ifx\SetFigFont\undefined%
\gdef\SetFigFont#1#2#3#4#5{%
  \reset@font\fontsize{#1}{#2pt}%
  \fontfamily{#3}\fontseries{#4}\fontshape{#5}%
  \selectfont}%
\fi\endgroup%
\begin{picture}(13330,4641)(2387,-5908)
\put(11010,-3741){\makebox(0,0)[b]{\smash{{\SetFigFont{17}{20.4}{\rmdefault}{\mddefault}{\updefault}$x$}}}}
\put(8710,-1491){\makebox(0,0)[b]{\smash{{\SetFigFont{17}{20.4}{\rmdefault}{\mddefault}{\updefault}$y$}}}}
\put(8701,-5836){\makebox(0,0)[b]{\smash{{\SetFigFont{20}{24.0}{\rmdefault}{\mddefault}{\updefault}$\lambda=0$}}}}
\put(15510,-3733){\makebox(0,0)[b]{\smash{{\SetFigFont{17}{20.4}{\rmdefault}{\mddefault}{\updefault}$x$}}}}
\put(13210,-1483){\makebox(0,0)[b]{\smash{{\SetFigFont{17}{20.4}{\rmdefault}{\mddefault}{\updefault}$y$}}}}
\put(13201,-5828){\makebox(0,0)[b]{\smash{{\SetFigFont{20}{24.0}{\rmdefault}{\mddefault}{\updefault}$\lambda>0$}}}}
\put(6508,-3737){\makebox(0,0)[b]{\smash{{\SetFigFont{17}{20.4}{\rmdefault}{\mddefault}{\updefault}$x$}}}}
\put(4208,-1487){\makebox(0,0)[b]{\smash{{\SetFigFont{17}{20.4}{\rmdefault}{\mddefault}{\updefault}$y$}}}}
\put(4199,-5832){\makebox(0,0)[b]{\smash{{\SetFigFont{20}{24.0}{\rmdefault}{\mddefault}{\updefault}$\lambda<0$}}}}
\end{picture}%

%% file: fig-vertex-3.pstex_t
\begin{picture}(0,0)%
\includegraphics{fig-vertex-3.pstex}%
\end{picture}%
\setlength{\unitlength}{3947sp}%
\begingroup\makeatletter\ifx\SetFigFont\undefined%
\gdef\SetFigFont#1#2#3#4#5{%
  \reset@font\fontsize{#1}{#2pt}%
  \fontfamily{#3}\fontseries{#4}\fontshape{#5}%
  \selectfont}%
\fi\endgroup%
\begin{picture}(13328,4634)(2389,-5909)
\put(11010,-3741){\makebox(0,0)[b]{\smash{{\SetFigFont{17}{20.4}{\rmdefault}{\mddefault}{\updefault}$x$}}}}
\put(8710,-1491){\makebox(0,0)[b]{\smash{{\SetFigFont{17}{20.4}{\rmdefault}{\mddefault}{\updefault}$y$}}}}
\put(8701,-5836){\makebox(0,0)[b]{\smash{{\SetFigFont{20}{24.0}{\rmdefault}{\mddefault}{\updefault}$\lambda=0$}}}}
\put(15510,-3742){\makebox(0,0)[b]{\smash{{\SetFigFont{17}{20.4}{\rmdefault}{\mddefault}{\updefault}$x$}}}}
\put(13210,-1492){\makebox(0,0)[b]{\smash{{\SetFigFont{17}{20.4}{\rmdefault}{\mddefault}{\updefault}$y$}}}}
\put(13201,-5837){\makebox(0,0)[b]{\smash{{\SetFigFont{20}{24.0}{\rmdefault}{\mddefault}{\updefault}$\lambda>0$}}}}
\put(6510,-3741){\makebox(0,0)[b]{\smash{{\SetFigFont{17}{20.4}{\rmdefault}{\mddefault}{\updefault}$x$}}}}
\put(4210,-1491){\makebox(0,0)[b]{\smash{{\SetFigFont{17}{20.4}{\rmdefault}{\mddefault}{\updefault}$y$}}}}
\put(4201,-5836){\makebox(0,0)[b]{\smash{{\SetFigFont{20}{24.0}{\rmdefault}{\mddefault}{\updefault}$\lambda<0$}}}}
\end{picture}%